\documentclass[a4paper,12pt,reqno]{amsart}

\usepackage[letterpaper, hmargin=1in, top=1in, bottom=1.2in, footskip=0.7in]{geometry}
\usepackage{cancel} 

\usepackage{rotating,pdflscape,xcolor,amssymb,subfigure,psfrag,
amsmath,eufrak,bbm,epsfig,amsthm,mathtools,fancyhdr,graphicx}
\definecolor{vdarkred}{rgb}{0.6,0,0.2}
\definecolor{vdarkblue}{rgb}{0,0.2,0.6}
\usepackage[pdftex, colorlinks, linkcolor=vdarkblue,citecolor=vdarkred]{hyperref}
\usepackage{a4wide,amscd}
\usepackage{tikz} 
 \usepackage[usenames,dvipsnames]{pstricks}
\usepackage[small]{caption}
\usepackage[normalem]{ulem} 

\newcommand{\me}{\mathrm{e}}
\newcommand{\trans}{\mathsf{T}}
\newcommand{\ld}{\ensuremath{\ldots}}

\newcommand{\sfu}{\mathsf{u}}
\newcommand{\sfU}{\mathsf{U}}
\newcommand{\sfv}{\mathsf{v}}

\newcommand{\sfy}{\mathsf{y}}
\newcommand{\sfz}{\mathsf{z}}
\newcommand{\sfX}{\mathsf{X}}

\newcommand{\cS}{\mathcal{S}}

\newcommand{\cX}{\mathcal{X}}

\newcommand{\bI}{\mathbf{I}}

\newcommand{\bi}{\mathbf{i}}\newcommand{\bj}{\mathbf{j}}

\newcommand{\Ga}{\Gamma} 
 
\newcommand{\ka}{\kappa} 

\newcommand{\si}{\sigma}
\newcommand{\lam}{\lambda}

\newcommand{\al}{\alpha}

\newcommand{\bet}{\beta}
\newcommand{\del}{\delta}

\newcommand{\cdote}{\,\cdot\,}

\newcommand{\beg}{\begin}
\newcommand{\en}{\end}

\newcommand{\bgt}{\begin{itemize}}
\newcommand{\ent}{\end{itemize}}
\newcommand{\ite}{\item}
\newcommand{\eqre}{\eqref}
\newcommand{\re}{\ref}
\newcommand{\la}{\label}

\newcommand{\brem}{\begin{rmk}}
\newcommand{\erem}{\end{rmk}}
\newcommand{\blem}{\begin{lem}}
\newcommand{\elem}{\end{lem}}
\newcommand{\bcor}{\begin{cor}}
\newcommand{\ecor}{\end{cor}}
\newcommand{\bTh}{\begin{Th}}
\newcommand{\eTh}{\end{Th}}
\newcommand{\bpropo}{\begin{propo}}
\newcommand{\epropo}{\end{propo}}

\newcommand{\op}{\operatorname} 

\newcommand{\Var}{\operatorname{Var}}

\newcommand{\Cov}{\operatorname{Cov}}

\newcommand{\Diag}{\operatorname{Diag}}

\newcommand{\Tr}{\operatorname{Tr}}

\newcommand{\ud}{\mathrm{d}}

\newcommand{\ninf}{\underset{n\to\infty}{\longrightarrow}}

\newcommand{\E}{\op{\mathbb{E}}}

\newcommand{\R}{\mathbb{R}}
\newcommand{\C}{\mathbb{C}}

\newcommand{\p}{\mathbb{P}}

\newcommand{\pro}{probability }

\newcommand{\f}{\frac}
\newcommand{\ff}{\frac{1}}
\newcommand{\lf}{\left}
\newcommand{\ri}{\right}

\newcommand{\st}{such that }

\newcommand{\ti}{\times}

\newcommand{\ste}{\, ;\, }

\newcommand{\eps}{\varepsilon}

\newcommand{\bbm}{\begin{bmatrix}}
\newcommand{\ebm}{\end{bmatrix}}
\newcommand{\bes}{\begin{equation*}}
\newcommand{\ees}{\end{equation*}}
\newcommand{\be}{\begin{equation}}
\newcommand{\ee}{\end{equation}}
\newcommand{\beqy}{\begin{eqnarray}}
\newcommand{\eeqy}{\end{eqnarray}}
\newcommand{\beq}{\begin{eqnarray*}}
\newcommand{\eeq}{\end{eqnarray*}}
\newcommand{\one}{\mathbbm{1}}

\newcommand{\ie}{i.e. }

\newcommand{\bpm}{\begin{pmatrix}}
\newcommand{\epm}{\end{pmatrix}}
\newcommand{\sbst}{\substack}

\newcommand{\cd}{\cdots}

\newcommand{\bpr}{\beg{proof}}
\newcommand{\epr}{\en{proof}}

\newcommand{\begenum}{\begin{enumerate}}
\newcommand{\enenum}{\end{enumerate}}
\newcommand{\MP}{\mathbb{MP}}

\newcommand{\NRM}[1]{{{\left\| #1\right\|}}} 
\newcommand{\Frobnorm}[1]{\left\| #1\right\|_{\op{F}}} 

\newcommand{\nsqt}{[n]^d_{<}}

\newcommand{\COV}{\mathcal{COV}}
\newcommand{\COR}{\mathcal{COR}}
\newcommand{\COVphi}{\COV^{\phi}}
\newcommand{\CORphi}{\COR^{\phi}}
\newcommand{\COVphika}{\COV^{\phi_\ka}}
\newcommand{\CORphika}{\COR^{\phi_\ka}}
\newcommand{\ovlphi}{\overline{\phi}}
\newcommand{\ntoinfty}{n\to\infty}
\newcommand{\Aij}{A_{\bi}}
\newcommand{\barUij}{\overline{U}_{\bi}}
\newcommand{\Mijone}{M^{(1)}_{\bi}}
\newcommand{\Mijtwo}{M^{(2)}_{\bi}}
\newcommand{\Mijthree}{M^{(3)}_{\bi}}
\newcommand{\Mijm}{M^{(l)}_{\bi}}
\newcommand{\Mklm}{M^{(l)}_{\bj}}
\newcommand{\barAij}{\overline{A}_{\bi}}
\newcommand{\xibullet}{\sfX_{i_1}(\bullet)}
\newcommand{\xjbullet}{\sfX_{i_d}(\bullet)}
\newcommand{\Do}{\mathbf{D}_0}
\newcommand{\ijinnsqt}{{\bi \in \nsqt}}
\newcommand{\klinnsqt}{{\bj \in \nsqt}}

\newtheorem{Th}{Theorem}[section]

\newtheorem{propo}[Th]{Proposition}

\newtheorem{lem}[Th]{Lemma}

\newtheorem{cor}[Th]{Corollary}
\theoremstyle{definition}
\newtheorem{rmk}[Th]{Remark}

\long\def\symbolfootnote[#1]#2{\begingroup
\def\thefootnote{\fnsymbol{footnote}}\footnote[#1]{#2}\endgroup}

\usepackage{listings}
\author[F. Benaych-Georges]{Florent Benaych-Georges}
\address{Florent Benaych-Georges: Capital Fund Management, 23 rue de l'Universit\'e, 75007 Paris, France}
\email{florent.benaych@gmail.com}

\author[T. Espana]{Tomas Espana}
\address{Tomas Espana: Chair of Econophysics and Complex Systems, Ecole Polytechnique, 91128 Palaiseau Cedex, France and LadHyX UMR CNRS 7646, \'Ecole Polytechnique, 91128 Palaiseau Cedex, France and Department of Operations Research and Financial Engineering, Princeton University, Princeton, NJ, 08540, USA}
\email{tomas.espana@princeton.edu}

\date{\today}

\title{Random Matrices and U-Statistics}
\begin{document}
\maketitle
\begin{abstract}
We introduce a family of coefficients based on U-statistics that generalize the notion of correlation and explore their properties in the large dimensional multivariate case, showing that in the null case of uncorrelated variables,  the spectrum of generalized correlation matrices is  distributed according to an affine transformation of the Marčenko-Pastur law.
\end{abstract}


 \section{Introduction}

Pearson's $\rho$, and Kendall's $\tau$ are two of the most commonly used metrics for quantifying the strength of the relationship between  random variables. Pearson's $\rho$ assesses the \textit{linear} relationships, making it particularly suitable for Gaussian data. In contrast, Kendall's $\tau$ is a \textit{rank}-based measure that offers greater robustness against heavy-tailed distributions. Understanding the unique characteristics of each coefficient is important for selecting the appropriate metric based on the application \cite{embrechts}. In this paper, however, we adopt a broader perspective by introducing a family of coefficients (see \eqre{defCOVphi}) that generalizes the notion of correlation (it encompasses both\footnote{Although closely related, Spearman’s rank correlation is not included in our generalization (see \cite{hoeffding}).} Pearson's $\rho$ and Kendall's $\tau$) and we explore their properties in the multivariate case, more specifically in the large dimensional regime for the null case of independent variables.  This motivates potential applications in data analysis techniques such as feature extraction (see, e.g., \cite{kendall_marko} for applications to portfolio construction in finance).

Random Matrix Theory (RMT) offers powerful tools for analyzing the properties of correlation matrices, especially in high-dimensional settings.   The  \emph{empirical spectral distribution} of a Hermitian matrix is the uniform distribution on its  eigenvalues.
  One of the purposes of RMT is to understand the asymptotic behaviour, as the dimension grows, of such distributions. In our case, we are interested in the asymptotic behaviour of these distributions when the entries of the matrix are  empirical  correlation coefficients (of the new, general type  introduced here at \eqre{defCOVphi}) of uncorrelated variables.

The first significant result in this area appeared with the analysis of the empirical spectral distribution of large sample covariance matrices in the null case. More specifically, the seminal work of \cite{marcenko_pastur} proved that the empirical spectral distribution of the sample covariance matrix of a white noise converges to a Marčenko-Pastur (MP) law when both the data dimension and the sample size grow (at the same rate). This is striking given the true covariance matrix of the white noise is the identity matrix, with spectrum concentrated at 1.  For Kendall's $\tau$, \cite{bandeirakendall} proved the convergence to an affine transformation of an MP law, elegantly using a variant of the Hoeffding decomposition.

Building on the work of \cite{bandeirakendall}, we show that the spectrum of generalized correlation matrices, with entries being U-statistics, are also distributed according to an affine transformation of a  MP law. This novel result, illustrated here by Figure \re{fig:is_c3},  highlights the intricate relationship between   MP distributions and sample correlation matrices of uncorrelated variables.

\begin{figure}[h!]
    \centering
    \includegraphics[width=0.55\textwidth]{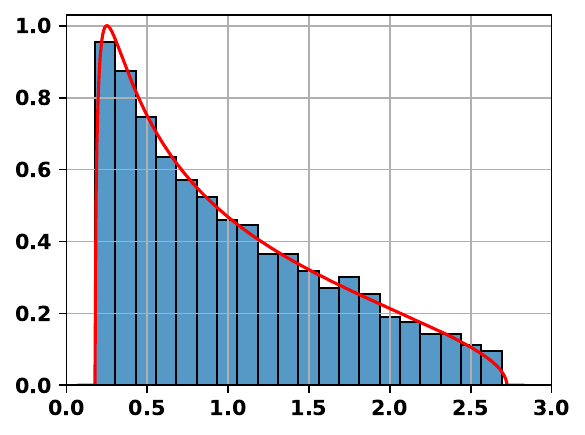}
    \caption{Histogram of the eigenvalues of $\CORphi$ when $d\!=\!2, \phi\!=\!\tanh, p\!=\!500$ and $q\!=\!0.5$, with i.i.d.\ standard Gaussian data. The superimposed red line is the density function of the Marčenko-Pastur law  $\MP_q^{t, 1-t}$ with $t\approx0.9$ (defined in Section \re{subsec:defMP}), as predicted by our Theorem \ref{maintheorem}}. 
    \label{fig:is_c3}
\end{figure}
 
\section{Main result}

For  $d\ge 2$ and $n$ positive integers, $\cX$ an arbitrary set,  $\phi:\cX^d\to\R$ \emph{antisymmetric}\footnote{A function $\phi:\cX^d\to\R$ is said to be \emph{antisymmetric} if for any $(x_1, \ldots, x_d)\in \cX^d$ and $\pi$ permutation of $[d]$, we have $\phi(x_{\pi(1)}, \ldots, x_{\pi(d)})=\op{sgn}(\pi)\phi(x_1, \ldots, x_d)$} and $\sfu, \sfv\in \cX^n$, 
we define the $\phi$-\emph{covariance}  $\COVphi(\sfu, \sfv)$   of $\sfu$ and $\sfv$  as \be\la{defCOVphi}
 \COVphi(\sfu, \sfv)=\ff{\binom{n}{d}}\sum_{(i_1, \ld, i_d)\in \nsqt} \phi(\sfu_{i_1}, \ldots, \sfu_{i_d}) \phi(\sfv_{i_1}, \ldots, \sfv_{i_d}),\ee 
where $[n]=\{1, \ld, n\}$ and 
 $$\nsqt=\{(i_1,\dots, i_d)\in [n]^d\ste i_1<\cdots<i_d\}.$$
When the denominator below is nonzero,   the $\phi$-\emph{correlation}  $\CORphi(\sfu,\sfv)$ of $\sfu$ and $\sfv$ as 
  \be\la{defCORphi} \CORphi(\sfu,\sfv)=\f{ \COVphi(\sfu,\sfv)}{\sqrt{ \COVphi(\sfu,\sfu) \COVphi(\sfv,\sfv)}}.\ee

Just as the covariance of two vectors extends to the notion of covariance matrix of a family of vectors, this notion extends to families of vectors: for $\sfu(1), \ld, \sfu(p)\in \cX^n$, we define the $p\ti p$ matrix
  \be\la{defMCOVphi}
   \COVphi(\sfu(1), \ld, \sfu(p))=\bbm  \COVphi(\sfu(k),\sfu(l))\ebm_{1\le k, l\le p}
   \ee
   and, when all involved denominators are nonzero, 
     \be\la{defMCORphi}
   \CORphi(\sfu(1), \ld, \sfu(p))=\bbm  \CORphi(\sfu(k),\sfu(l))\ebm_{1\le k, l\le p}
   \ee

   The following theorem, illustrated by Figure \re{fig:is_c3}, is the main result of this text.   
   
   \bTh\la{maintheorem} Let $d\ge 2$ be an integer. For each $n$, let $p$ be a positive integer, $$\lf(\sfX_i(k)\ri)_{\sbst{1\le i \le n\\ 1\le k\le p}}$$ be an array of i.i.d. random variables taking values in a measured space,
and $\phi:\cS^d\to\R$ be an antisymmetric measurable function (all, except $d$,  implicitly depending on $n$). We introduce  \be\la{eq:introxxprime}\cX_1, \ldots, \cX_d\ee some independent copies of $\sfX_1(1)$ and the function 
   \be\la{eq:defovlphi}\ovlphi(x)=\E\lf[\phi(x, \cX_2, \ld, \cX_d)\ri] \ee
    (so that the random variables $\cX_1, \ldots \cX_d$ and the function $\ovlphi$ implicitly depend on $n$).
   Suppose that   \bgt\ite the ratio $p/n$ tends to a finite limit $q\ge 0$ as $\ntoinfty$, 
   \ite   $\E\lf[\phi(\cX_1, \ldots, \cX_d)^{2d+4}\ri]$
   is bounded by a finite constant   independent of $n$, 
   \ite  the limits  $$\al=\lim_{\ntoinfty}d\E\lf[\,\ovlphi(\cX_1)^2\, \ri],\qquad \beta =\lim_{\ntoinfty}\E\lf[\phi(\cX_1, \ld, \cX_d)^2\ri]$$  exist and $\bet>0$.   \ent
   Then, as $\ntoinfty$, for $\sfX(j)=(\sfX_i(j))_{1\le i\le n}$, the empirical spectral distribution\footnote{The \emph{empirical spectral distribution} of a Hermitian matrix is the uniform \pro measure on its  eigenvalues (counted with multiplicity).}  of 
     $$\COVphi(\sfX(1), \ld, \sfX(p))$$ converges weakly in probability to the \pro measure\footnote{The \pro distributions $\MP_q^{\al, \bet}$ are defined in Section \re{subsec:defMP}. These are afffine transformations of Marčenko-Pastur laws.} $\MP_q^{\al, \bet-\al}$  and the empirical spectral distribution of 
     $$\CORphi(\sfX(1), \ld, \sfX(p))$$ converges weakly in probability to  $\MP_q^{t, 1-t}$ for $t=\f\al\bet$.
   \eTh

   \beg{rmk} Note that as dimension $p$,  the space $\cS$, $\phi$, the distribution of $\sfX_i(k)$ might depend on $n$. The dependence in $n$ is dropped from notations for brevity. 
The only parameters that do not depend on $n$ are the fixed integer $d$, the limits $q$, $\al$ and $\bet$, and the implicit uniform bound  on $\E\lf[\phi(\cX_1, \ld, \cX_d)^{2d+4}\ri]$.
   \en{rmk}

   \beg{rmk} Note that the hypothesis that the limit of $\E\lf[\phi(\cX_1, \ld, \cX_d)^2\ri]$ is positive implies that with probability tending to one, none of the $\COVphi(\sfX(k), \sfX(k))$ (for $1\le k\le p$) is null (this can be seen as a consequence of Lemmas \re{cutofflemma} and \re{lemconcdiagmat}), so that    $$\CORphi(\sfX(1), \ld, \sfX(p))$$ is actually well defined.
   \en{rmk}

  \beg{rmk}When $d\ge 3$, by antisymmetry of $\phi$ and exchangeability of the $\cX_i$'s, the function $\ovlphi$ defined at \eqre{eq:defovlphi} is null.  Hence, in this  case, the statements of the theorem come down  to say that the spectrum of $\COVphi(\sfX(1), \ld, \sfX(p))$ (resp. of $\CORphi(\sfX(1), \ld, \sfX(p))$) concentrates asymptotically in any neighborhood of $\bet$ (resp.\ of $1$), which is a non-trivial result. This opens the way to further investigations on the use of values of $d$ that are larger than $2$ to identify the null case of uncorrelated variables  via the spectrum of such matrices, even at high ratios $p/n$.
  \en{rmk}
   
      \beg{rmk} This generalization of correlation coefficients was inspired by \cite{daniels}, when $d=2$. In this case, a classical example of antisymmetric function is $\phi(x,y)=\psi(y-x)$ for $\psi:\R^d\to\R$ odd. Note that when $\psi(x) = x$, $\frac{1}{2}\COV^\phi(u, v)$ is exactly the standard unbiased estimator of covariance so that $\CORphi(u,v)$ is the Pearson's $\rho$ correlation coefficient. Also, when $\psi(x) = \text{sgn}(x)$, $\CORphi(u,v)$ corresponds to the Kendall's $\tau$ correlation coefficient. Thus, Theorem~\ref{maintheorem} recovers both the classical Marčenko-Pastur theorem and the result proved in \cite{bandeirakendall} (with a relaxation of the hypotheses). For  general $\cS$, lots of examples for $d\ge 3$  can stem from the remark that  for any antisymmetric function $\phi_0:\cS^2\to\R$, the function $\phi: \cS^d\to \R$ defined by $$\phi(x_1, \ld, x_d)=\prod_{(i,j)\in [d]_<^2} \phi_0(x_i, x_j)$$ is antisymmetric.
   \en{rmk}   
   
   \beg{rmk}The cornerstone of the proof of the theorem is   Lemma \re{mainlemma1}, whose idea has been inspired by  \cite{bandeirakendall}. It can be summed-up as follows: in  each entry \begin{equation*}  \qquad\qquad \ff{\binom{n}{d}}\sum_{(i_1, \ld, i_d) \in\nsqt} \phi(\sfX_{i_1}(k), \ld, \sfX_{i_d}(k)) \phi(\sfX_{i_1}(l), \ld,\sfX_{i_d}(l))    \qquad\qquad \text{($1\le k,l\le p$)}\end{equation*} of the matrix $\COVphi(\sfX(1), \ld, \sfX(p))$,
  one can replace the average over $\{i_2, \ld, i_d\}$ by the corresponding theoretical expectation without changing the spectrum of the matrix at the macroscopic level. This is surprising given the whole Marčenko-Pastur theory relies on the fact that in an empirical covariance matrix, one cannot replace empirical means by their theoretical values without triggering a dramatic macroscopic change in the spectrum.  \en{rmk}

\beg{rmk}
  For any $\sfu, \sfv\in \cX^n$ and $d\geq 1$,  the $\phi$-\emph{covariance}  $\COVphi(\sfu, \sfv)$ is in fact a U-statistic of degree $d$ with (symmetric) kernel $h(\sfz_1,\ldots,\sfz_d) = \phi(\sfu_1,\ldots,\sfu_d)\phi(\sfv_1,\ldots,\sfv_d)$ where $\sfz_i = (\sfu_i, \sfv_i)$. This property infers to $\phi$-\emph{covariances} many properties like law of large numbers, asymptotic normality, law of iterated logarithms and large deviations \cite{lee2019u, de2012decoupling}. \en{rmk}
   
   \section{Proofs}
   \subsection{Notation} Let $M^\trans$ denote the transpose of a matrix $M$. 
   For $u, v\in \R^p$, we set  $u\otimes v=uv^\trans$, where $u, v$ are considered as column vectors (so that  $u\otimes v$ is a  $p\ti p$ matrix). 
   Let us also introduce the operator $\Do$ that maps any square matrix to the same matrix, where all diagonal entries have been set to zero and $\Diag$, the operator that maps any vector to the square diagonal matrix with diagonal entries the coordinates of the vector. We also use the notation $\Diag(M)=M-\Do(M)$ for $M$ a square matrix.
    At last, we define the Frobenius norm $\Frobnorm{M}$ of a real matrix $M$ as follows \be\la{deffrobnorm} \Frobnorm{M}=\sqrt{\Tr MM^\trans}.\ee
   
         \subsection{Proof of Theorem \re{maintheorem}}Let us prove Theorem \re{maintheorem}
         based on Lemmas \re{cutofflemma},  \re{mainlemma1}, \re{mainlemma2}, and \re{lemconcdiagmat},  that will be proved below. 
         
         Remember that $\phi$  depends implicitly on $n$. For $\ka >0$, we define $\phi_\ka:\cS^d\to\R$ as follows: $$\phi_\ka(x_1, \ld, x_d)=\op{sign}(\phi(x_1, \ld, x_d))\min(|\phi(x_1, \ld, x_d)|, n^\ka).$$  Then, $\ovlphi_\ka$ is defined as $\ovlphi$: $\ovlphi_\ka(x)=\E\lf[\phi_\ka(x, \cX_2, \ld, \cX_d)\ri]$.
      Note that $\phi_\ka$ is bounded by $n^\ka$ and $\phi(x_1, \ld, x_d)=\phi_\ka(x_1, \ld, x_d)$ when $|\phi(x_1, \ld, x_d)|\le n^\ka$. 
      \blem\la{cutofflemma}
      There is an exponent $\ka\in (0, 1/2)$ \st
      \be\la{probatendstozero} \p\lf(\exists (i_1, \ld, i_d, k)\in [n]^d\ti [p] \text{ such that } |\phi(\sfX_{i_1}(k), \ld, \sfX_{i_d}(k))|>n^\ka\ri)\ninf 0
      \ee
      and for any fixed integer $k\ge 3$,       \be\la{notoolargemoments}
      \f{\E\lf[ |\ovlphi_\ka(\cX_1)|^k\ri]}{n^{\f{k}2-1}}\ninf 0.
      \ee
      \elem

      By \eqre{probatendstozero}, the \pro that $$\COVphi(\sfX(1), \ld, \sfX(p))\ne \COVphika(\sfX(1), \ld, \sfX(p))$$ or $$\CORphi(\sfX(1), \ld, \sfX(p))\ne \CORphika(\sfX(1), \ld, \sfX(p))$$ tends to zero as $\ntoinfty$, so given we prove a convergence in probability, we can switch from $\phi$ to $\phi_\ka$, \ie, suppose that for a certain $\ka\in (0, 1/2)$,  \be\la{boundphiaftercutoff}
    |\phi|\le n^\ka  \ee and for any $k\ge 3$,   \be\la{notoolargemomentscutoff}
      \f{\E\lf[ |\ovlphi(\cX_1)|^k\ri]}{n^{\f{k}2-1}}\ninf 0.
      \ee
      
      For $i\in [n]$, we define the vectors  of $\R^p$ \beqy
   \nonumber\sfX_i(\bullet)&=&(\sfX_i(1), \ld, \sfX_i(p)) \\
\la{eq:defUi} \sfU_i&=&\ovlphi(\sfX_i(\bullet)) 
\eeqy where, in   \eqre{eq:defUi},   $\ovlphi$ is  applied entry-wise to each of the $p$ components.

  The following equations are obvious from the fact that $\phi$ is antisymmetric  and the $\lf(\sfX_i(k)\ri)_{\sbst{1\le i \le n\\ 1\le k\le p}}$ are i.i.d.\, and will be used extensively in the following: for any $\ijinnsqt$, and any $k\in [p]$, 
$$\E[\phi(\sfX_{i_1}(k), \ld, \sfX_{i_d}(k))]=\E[\ovlphi(\sfX_{i_2}(k))]=\cdots=\E[\ovlphi(\sfX_{i_d}(k))]=0.$$

\blem\la{mainlemma1} Let us define the $p\times p$ matrix \be\la{defGa}\Ga=\f dn\sum_{i=1}^n \sfU_i\otimes \sfU_i+\lf(\E\lf[\phi(\cX_1, \ld, \cX_d)^2\ri]-d\E\lf[\,\ovlphi(\cX_1)^2\, \ri]\ri)\bI_p.\ee
Then $$\E \lf[\Frobnorm{\COVphi(\sfX(1), \ld, \sfX(p))-\Ga}^2\ri]=O(1).$$
\elem   

\blem\la{mainlemma2} The empirical spectral law of the  $p\times p$ matrix   $\Ga$ from  \eqre{defGa} converges weakly in \pro to $\MP_q^{\al, \bet-\al}$.  Besides, we have $\E[\Frobnorm{\Ga}^2]=O(p).$
\elem

Now, the proof of the convergence of the empirical spectral law of 
     $$\COVphi(\sfX(1), \ld, \sfX(p))$$  to $\MP_q^{\al, \bet-\al}$ is a direct application   of Lemmas \re{mainlemma1}, \re{mainlemma2} and \re{lemFrobeniusapproximation}.
     
     \blem\la{lemconcdiagmat} For $\COVphi(\cdote, \cdote)$ defined at \eqre{defCORphi}, we have, for any $k\in [p]$ and $r\ge 0$,  $$\p\lf(\lf|\COVphi(\sfX(k), \sfX(k))-\E\lf[\phi(\cX_1, \ld, \cX_d)^2\ri]\ri|\ge  r\ri)\le 2\me^{-2r^2/\si^2}$$ for    $\si^2=d^{2d+2}n^{2\ka -1}$.
     \elem
     
     Let us now prove the convergence of the empirical spectral law of 
     $$\CORphi(\sfX(1), \ld, \sfX(p))$$  to $\MP_q^{t, 1-t}$ for $t=\f\al\bet$. Note that $$\CORphi(\sfX(1), \ld, \sfX(p))=D^{-1/2}\COVphi(\sfX(1), \ld, \sfX(p))D^{-1/2}$$ for $D$ the diagonal matrix defined by  $$D=\bpm \COVphi(\sfX(1), \sfX(1))&&\\ &\ddots&\\ && \COVphi(\sfX(p), \sfX(p))\epm$$
     By Lemma \re{lemconcdiagmat}, there is a constant $C$ depending only on $d$ \st  for any $k\in [p]$,  $$\p\lf(\lf|\COVphi(\sfX(k), \sfX(k))-\E\lf[\phi(\cX_1, \ld, \cX_d)^2\ri]\ri|\ge  C n^{\ka/2-1/4}\ri)\le 2\me^{-n^{1/2-\ka}},$$ hence by the union bound,  
     $$\p\lf(\NRM{D- \E\lf[\phi(\cX_1, \ld, \cX_d)^2\ri]\bI_p}\ge  C n^{\ka/2-1/4}\ri)
     \le 2p\me^{-n^{1/2-\ka}},$$
so that the hypotheses of Lemma \re{lemFrobeniusapproximationwithscalarmatrix} are satisfied for $M_n=     \COVphi(\sfX(1), \ld, \sfX(p))$, $\lam=\bet$ and $f(x)=\one_{x>0}x^{-1/2}$, which concludes the proof. 

    \subsection{Proof of  Lemma \ref{cutofflemma}}
      By hypothesis, $\E\lf[\phi(\cX_1, \ld, \cX_d)^{2d+4}\ri]$
   is bounded by a finite constant independent of $n$, so there is a constant $C$ \st uniformly in $\ka>0$,  $n$ and $(i_1, \ld,i_d,k)\in [n]^d \ti[p]$,  $$\p\lf( |\phi(\sfX_{i_1}(k), \ld, \sfX_{i_d}(k))|>n^\ka\ri)\le   C n^{-(2d+4)\ka}$$
   Hence by the union bound,  \eqre{probatendstozero} is satisfied as soon as 
   \be\la{cond1onkappa}d+1-(2d+4)\ka <0.\ee
   Note that by the hypothesis on $\E\lf[\phi(\cX_1, \ld, \cX_d)^{2d+4}\ri]$ again, \eqre{notoolargemoments} is always satisfied for $k\in\{3, \ld,  2d+4\}$. For $k\ge  2d+5$, we have $$\E\lf[ |\ovlphi_\ka(\cX_1)|^k\ri]\le n^{k\ka },$$ so that \eqre{notoolargemoments}  is satisfied as soon as $k\ka <\f{k}2-1$, \ie 
   \be\la{cond2onkappa}\ka <\ff2-\ff{k},\ee which is satisfied for any $k\ge 2d+5$ as soon as it is satisfied for $k=2d+5$. 
   Joining conditions \eqre{cond1onkappa} and \eqre{cond2onkappa}, we get that the lemma is proved if we can find $\ka \in (0, 1/2)$ \st $$\f{d+1}{2d+4}<\ka<\ff2-\ff{2d+5}.$$ Given $$\f{d+1}{2d+4}+\ff{2d+5}<\f{d+1}{2d+4}+\ff{2d+4}=\ff2,$$ this concludes the proof.

   \subsection{Proof of  Lemma \ref{mainlemma1}}\leavevmode\par
   
   For   $\bi = (i_1, \ld, i_d)\in \nsqt$, we define the vector  of $\R^p$ \beq
\la{eq:defAij}\Aij&=&\phi(\xibullet, \ld, \xjbullet) ,\eeq where   $\phi$ is  applied entry-wise to each of the $p$ components.

Note   that, by the very definition  \eqre{defCOVphi} of  $\COVphi$, we have  $$\COVphi(\sfX(1), \ld, \sfX(p))=\ff{\binom{n}{d}}\sum_\ijinnsqt\Aij\otimes \Aij$$
Note also the following relations  between $\Aij$'s conditional expectations and the $\sfU_{i}$'s defined at \eqre{eq:defUi}: \begin{equation}\label{2408250837}
\begin{aligned}
\E\!\left[\Aij \mid \sfX_{i_1}(\bullet)\right] &= \sfU_{i_1},\\
\E\!\left[\Aij \mid \sfX_{i_l}(\bullet)\right] &= -\,\sfU_{i_l}, \qquad \forall\, l\in\{2,\ldots,d\}.
\end{aligned}
\end{equation}

Then,  for   $\bi = (i_1, \ld, i_d)\in \nsqt$, we define the vector  of $\R^p$ \beqy
\la{eq:defAbarij}\barAij&=& \Aij-\lf(\sfU_{i_1}\!-\!\sfU_{i_2}\!-\!\dotsb\!-\!\sfU_{i_d}\ri).
\eeqy Henceforth, we denote $\barUij = \sfU_{i_1}\!-\!\sfU_{i_2}\!-\!\dotsb\!-\!\sfU_{i_d}$.

Then, we introduce the matrices $\Mijone, \Mijtwo, \Mijthree$ defined as follows
\begin{align}
\Mijone
  &= \Diag\!\lf(\phi^2(\xibullet,\ld,\xjbullet)\ri)
     + \Do\lf(\barUij \otimes \barUij \ri) \la{Mijone} \\
\Mijtwo
  &= \Do\lf(\barAij \otimes \barUij \ri) \notag \\
\Mijthree
  &= \Do\lf(\barAij \otimes \barAij\ri) \notag
\end{align}where, in   \eqre{Mijone},   $\phi^2$ is  applied entry-wise  to each of the $p$ components.

   Then, using the identity, stemming from \eqre{eq:defAbarij}, that $\Aij=\barUij + \barAij$, we have the decomposition
   \beq 
   \Aij\otimes \Aij&=& \Diag(\phi^2(\xibullet, \ld, \xjbullet))+\Do(  \Aij\otimes \Aij)\\
   &=&\Diag(\phi^2(\xibullet, \ld, \xjbullet))+\Do\lf(\barUij \otimes \barUij \ri)+\\ &&\Do\lf(\barAij \otimes \barUij \ri)+\Do\lf(\barUij \otimes\barAij\ri)+\\ &&\Do\lf(\barAij \otimes\barAij  \ri)\\
   &=& \Mijone+ \Mijtwo+(\Mijtwo)^\trans+\Mijthree
   \eeq
   
   Then, the proof of Lemma \ref{mainlemma1} follows directly from Fact \re{Fact23082514h} and Fact \re{Fact23082514hM23} below, using the triangular inequality for the norm $\Frobnorm{\cdote}$  and Inequality \eqre{ineqsumsquares}.  
   \beg{fact}\la{Fact23082514h}For $\Ga$ as in \eqre{defGa}, we have $$\E\lf[\Frobnorm{
  \ff{\binom{n}{d}}\sum_\ijinnsqt\Mijone -\Ga}^2
   \ri]=O(1).$$
   \en{fact}

\beg{fact}\la{Fact23082514hM23}For any $l\in \{2, 3\}$, 
$$ \E\lf[\Frobnorm{\sum_\ijinnsqt\Mijm}^2\ri]=O(n^{2(d-1)}p^2)$$
\en{fact}

\emph{Proof of Fact \re{Fact23082514h}.} Using Lemma \re{polynomialidentity}, we have 
\beq \sum_\ijinnsqt \barUij \otimes \barUij &=&\binom{n\!-\!1}{d\!-\!1}\sum_{i=1}^n \sfU_i\otimes \sfU_i + \sum_{\sbst{1 \leq i,j\leq n \\ i \neq j}} \left(\binom{n\!-\!2}{d\!-\!2} - 2\binom{n\!-\!\min(i,j)\!-\!1}{d\!-\!2}\right) \sfU_i \otimes \sfU_j 
\eeq 
Hence,
\begin{align}
\ff{\binom{n}{d}}\sum_\ijinnsqt \Mijone - \Ga
&= \ff{\binom{n}{d}}\sum_\ijinnsqt
   \Big(\Diag(\phi^2(\xibullet,\ld,\xjbullet))
      - \E\lf[\phi(\cX_1,\ld,\cX_d)^2\ri]\bI_p\Big) \notag \\
&\quad - \f dn\sum_{i=1}^n \Big( \Diag(\sfU_i\otimes \sfU_i)
      - \E\lf[\,\ovlphi(\cX_1)^2\,\ri]\bI_p\Big) \notag \\
&\quad + \ff{\binom{n}{d}}\Do\, \sum_{\sbst{1 \leq i,j\leq n \\ i \neq j}} \left(\binom{n\!-\!2}{d\!-\!2} - 2\binom{n\!-\!\min(i,j)\!-\!1}{d\!-\!2}\right) \sfU_i \otimes \sfU_j  \la{decdiffsum2308251}
\end{align}

The first term of the RHT of \eqre{decdiffsum2308251} is treated as follows, using \eqre{keyupperbound2308252}:
\begin{align}
&\E\!\Big[\Big\|\ff{\binom{n}{d}}
   \sum_\ijinnsqt\!\Big(\Diag(\phi^2(\xibullet,\ld,\xjbullet))
   - \E\lf[\phi(\cX_1,\ld,\cX_d)^2\ri]\bI_p\Big)\Big\|_F^{\!2}\Big] = \notag\\
\la{keyupperbound2308252applic1}
&\sum_{k=1}^p \Var\!\left(
   \ff{\binom{n}{d}}
   \sum_\ijinnsqt  \phi^2(\sfX_{i_1}(k),\ld,\sfX_{i_d}(k))
   \right)
\;\le\; \frac{2pd^d}{n}\,\Var\!\big(\phi(\cX_1,\ld,\cX_d)^2\big),
\end{align} using the fact that $\binom{n}{d}\ge (n/d)^d$.

The second term of the RHT of \eqre{decdiffsum2308251} is treated as follows:
\beqy\nonumber &&
\E\lf[\Frobnorm{\f dn\sum_{i=1}^n \lf( \Diag(\sfU_i\otimes \sfU_i)-\E\lf[\,\ovlphi(\cX_1)^2\, \ri]\bI_p\ri)}^2\ri]=
\\ \la{23082512h43}&&\sum_{k=1}^p\Var\lf(\f dn\sum_{i=1}^n \ovlphi(\sfX_i(k))^2 \ri)\;=\; \f{d^2p}n \Var(\ovlphi(\cX_1)^2)
\eeqy 
The third term of the RHT of \eqre{decdiffsum2308251} is treated as follows:
\begin{align}
&\E\lf[\Frobnorm{\ff{\binom{n}{d}}\Do
\sum_{\sbst{1 \le i,j \le n \\ i \neq j}}
\left(\binom{n\!-\!2}{d\!-\!2} - 2\binom{n\!-\!\min(i,j)\!-\!1}{d\!-\!2}\ri)
\sfU_i \otimes \sfU_j }^2\right] \notag\\
&= \ff{\binom{n}{d}^2}
\sum_{\sbst{1\le k, l\le p\\ k\ne l}}
\E\Bigg(\Bigg(\sum_{\sbst{1\le i, j\le n\\ i\ne j}} \left(\binom{n\!-\!2}{d\!-\!2} - 2\binom{n\!-\!\min(i,j)\!-\!1}{d\!-\!2}\ri)
\ovlphi(\sfX_i(k))\ovlphi(\sfX_j(l))\Bigg)^2\Bigg) \notag\\
&\leq \frac{\binom{n\!-\!2}{d\!-\!2}^2}{\binom{n}{d}^2}
\sum_{\sbst{1\le k, l\le p\\ k\ne l}}
\E\Bigg(\sum_{\sbst{1\le i, j\le n\\ i\ne j}}
\ovlphi(\sfX_i(k))^2\,\ovlphi(\sfX_j(l))^2\Bigg) \notag\\
\la{23082512h43bis}
&\leq \frac{Kp^2}{n^2}\,\E\big[\ovlphi(\cX_1)^2\big]^2,
\end{align}
for a constant $K > 0$ depending only on $d$. Then, we conclude the proof of Fact \re{Fact23082514h} starting from \eqre{decdiffsum2308251},  using the triangular inequality for the norm $\Frobnorm{\cdote}$  and \eqre{ineqsumsquares}, and treating   the three terms on the RHT thanks to respectively \eqre{keyupperbound2308252applic1}, \eqre{23082512h43}, and  \eqre{23082512h43bis}.
\hfill $\square$\\

\emph{Proof of Fact \re{Fact23082514hM23}.} 
Let $l\in\{2,3\}$. First, for any $\ijinnsqt$ and any $\klinnsqt$, $$\bigl|\{i_1, \ld, i_d\} \cap \{j_1, \ld, j_d\}\bigr| \leq 1\implies\E\Tr (\Mijm(\Mklm)^\trans) =0$$ (this stems  from Tower property, \eqre{2408250837} and   the fact that the $\sfU_i$'s are centered, as $\phi$ is antisymmetric). 
It follows that for $$\mathcal{C}_d = 
\Bigl\{ (\mathbf{i},\mathbf{j}) \in \nsqt \times \nsqt :
\bigl|\{i_1,\dots,i_d\} \cap \{j_1,\dots,j_d\}\bigr| \ge 2 \Bigr\},$$  we have \beq \E\lf[\Frobnorm{\sum_\ijinnsqt\Mijm}^2\ri]&=&\sum_{(\bi, \bj) \in \mathcal{C}_d}\sum_{1\le k, l\le p}\E\lf[(\Mijm)_{k,l}(\Mklm)_{k,l}\ri].
\eeq
This yields the conclusion, given that $\bigl| \mathcal{C}_d\bigr| = O(n^{2(d-1)})$ and that the second moments of the entries of $\Mijm$ are dominated by the fourth moments of $\phi(\cX_1, \ldots, \cX_d)$, which are supposed to be bounded. 
\hfill $\square$\\

  \subsection{Proof of  Lemma \ref{mainlemma2}}When $q>0$, the convergence is a direct application of \cite[Th. 3.2]{FloThierryCabDuv}. Indeed, $$\ff{n}\sum_{i=1}^n \sfU_i\otimes \sfU_i=\ff{n}MM^\trans$$ for $$M=[\ovlphi(\sfX_i(k))]_{\substack{1\le i\le n\\ 1\le k\le p}}$$ and we have the bound \eqre{notoolargemomentscutoff} on the moments of the entries of $M$.  The fact that $\E[\Frobnorm{\Ga}^2]=O(p)$ is for example a consequence of  Proposition \re{flatMPcase} (used with \eqre{ineqsumsquares}). 
  
  For the case $q=0$, we simply use Proposition \re{flatMPcase} with Lemma \re{lemFrobeniusapproximation}.
  
    \subsection{Proof of  Lemma \ref{lemconcdiagmat}}The lemma is a direct application	of Lemma \re{McDiarmidLemma} for $$f(x_1, \ld, x_n)=\ff{\binom{n}{d}}\sum_\ijinnsqt \phi(x_{i_1}, \ldots, x_{i_d})^2    $$with the observation that given $|\phi|\le n^\ka$, the numbers $\si_i$ of \eqre{defsiisquare} satisfy 
    $$\si_i\le \f{dn^{d-1}n^\ka }{\binom{n}{d}} \le d^{d+1}n^{\ka-1} $$ (we used the fact that $\binom{n}{d}\ge (n/d)^d$), so that $$\sum_{i=1}^n \si_i^2\le d^{2d+2}n^{2\ka -1}.$$
  
 \section{Appendix}

\subsection{Marčenko-Pastur distributions $\MP_q$ and $\MP_q^{\al, \bet}$}\la{subsec:defMP}
 
 For $q\ge 0$, we define the \emph{Marčenko-Patsur distribution} $\MP_q$ with \emph{parameter} $q$ as follows:
 \be\la{eq:defMP}\MP_q=\beg{cases}\del_1&\text{ if }q=0,\\
  \f{\sqrt{(\lam_+-x)(x-\lam_-)}}{2\pi q x}\ud x+\one_{q>1}(1-q^{-1})\del_0&\text{ if }q>0,
 \en{cases}
 \ee
 for $\lam_\pm=(1\pm\sqrt{q})^2$.
 
 For $\al, \bet\in \R$, we set $\MP_q^{\al, \bet}$ to be the push-forward of $\MP_q$ by the map $x\mapsto \al x +\bet$.
 
 \subsection{A polynomial identity}
 The following identity is one of the keys of our proof. 
 \blem\la{polynomialidentity} Let  $ 2\le d\le n$. In the space $\C[X_1,Y_1, \ld, X_n, Y_n]$ of polynomials in the indeterminates $X_1,Y_1, \ld, X_n, Y_n$, the following identity holds:
\beq
&&\sum_\ijinnsqt (X_{i_1}\!-\!X_{i_2}\!-\!\cd\!-\!X_{i_d})(Y_{i_1}\!-\!Y_{i_2}\!-\!\cd-\!Y_{i_d})\\
&& \qquad \;=\;\binom{n\!-\!1}{d\!-\!1}\sum_{i=1}^n X_i Y_i + \sum_{\substack{1 \le i,j \le n\\ i \ne j}} \left(\binom{n\!-\!2}{d\!-\!2}\!-\!2\binom{n\!-\!\min(i,j)\!-\!1}{d\!-\!2}\right) X_i Y_j.
\eeq
\elem 
\bpr
When expanding the outer product, each diagonal term $X_iY_i$ appears for every $d$-subset containing $i$, i.e., $\binom{n-1}{d-1}$ times. This yields the first term of the RHT. For the second term, let $i\neq j$ and, w.l.o.g, $i<j$. The total number of $d$-subsets containing $\{i,j\}$ is $\binom{n-2}{d-2}$. Among these, those where $i$ is the smallest element of the subset consists in picking the other $d\!-\!2$ elements from the $n\!-\!i\!-\!1$ numbers in $\{i\!+\!1,\ldots,n\} \backslash \{j\}$, giving $\binom{n-i-1}{d-2}$ terms. For $i$ in a $d$-subset, we denote $\sf{rk}(i)$ its rank among the subset and consider \[
c_r = 
\begin{cases}
+1, & r = 1, \\
-1, & r = 2, \ldots, m.
\end{cases}
\]
Then, if $i$ is the smallest of the $d$-subset, $c_{\sf{rk}(i)}c_{\sf{rk}(j)}\!=\!-\!1$ and $c_{\sf{rk}(i)}c_{\sf{rk}(j)}\!=\!1$ otherwise. Therefore, for $i<j$, the term in front of $X_iY_j$ is,
\[ (+1) \cdot \left(\binom{n\!-\!2}{d\!-\!2}\!-\!\binom{n\!-\!i\!-\!1}{d\!-\!2}\right)+(-1)\cdot \binom{n\!-\!i\!-\!1}{d\!-\!2} = \binom{n\!-\!2}{d\!-\!2}\!-\!2\binom{n\!-\!i\!-\!1}{d\!-\!2}
\]
The result follows by symmetry for general $i\neq j$.
\epr

 \subsection{Approximation and Frobenius norm}
The following lemma follows directly from  \cite[Cor. A.41]{bai-silver-book}. 
\blem\la{lemFrobeniusapproximation}
For $n\ge 1$, let $M_n, M_n'$ be symmetric real matrices, both having size $p_n\times p_n$, with $p_n\to\infty$ as $\ntoinfty$. Suppose that the empirical spectral distribution of $M_n$ converges weakly in probability to a non random limit probability measure $\mu$ and that, for $\Frobnorm{\cdote}$ the Frobenius norm defined at \eqre{deffrobnorm}, $$\f{\Frobnorm{M_n'-M_n}^2}{p_n}.$$ converges in \pro to zero. 
 Then the empirical spectral distribution of $M_n'$ converges weakly in probability to $\mu$.
\elem

Let us denote by $\NRM{\cdote}$ the operator norm associated with the canonical Euclidian norm. In the lemma below, $f(D_n)$ has to be understood as the result of the functional calculus (\ie $f$ is applied to the eigenvalues of $D_n$ and the eigenvectors remain unchanged).  
 \blem\la{lemFrobeniusapproximationwithscalarmatrix}For $n\ge 1$, let $M_n, D_n$ be symmetric real matrices, both having size $p_n\times p_n$, with $p_n\to\infty$ as $\ntoinfty$.  Suppose that  
 \bgt\ite[(i)] the empirical spectral distribution of $M_n$ converges weakly in probability to a non random limit probability measure $\mu$,
 \ite[(ii)]  $\E\lf[\Frobnorm{M_n}^2\ri]=O(p_n)$,
 \ite[(iii)] there is $\lam \in \R$ \st $\NRM{D_n-\lam \bI_{p_n} }$ converges in \pro to zero as $\ntoinfty$.\ent
 Then for any measurable  $f:\R\to\R$ continuous at $\lam$,  the empirical spectral distribution of $f(D_n)M_nf(D_n)$ converges weakly in probability to the push-forward of $\mu$ by the function $x\mapsto f(\lam)^2 x$. 
 \elem
 
 \bpr Note that (iii) simply means that any neighborhood of $\lam$ will, with \pro tending to one as $\ntoinfty$, contain the whole spectrum of $D_n$. Hence by continuity of $f$ at $\lam$, we can suppose that $f$ is the identity function. Then,
  by Lemma \re{lemFrobeniusapproximation}, one has to prove that  $$\f{\Frobnorm{D_nM_nD_n-\lam^2M_n}^2}{p_n}.$$ converges in \pro to zero. 
   For $\eps_n=D_n-\lam \bI_{p_n}$, by \eqre{ineqsumsquares}, we have $$\Frobnorm{D_nM_nD_n-\lam^2M_n}^2\le 8|\lam|\Frobnorm{\eps_nM_n}^2+4\Frobnorm{\eps_n^2M_n }^2.$$Then, one concludes using the non-commutative H\"older  inequality  \cite[Appendix A.3]{alice-greg-ofer} and the fact that convergence in \pro to zero is robust to multiplication by $L^1$-bounded sequences.
 \epr
 
Now comes a  statement that we need for the case $q=0$ and that we did not find in literature. Note that this elementary computation shows, without the need of the whole Marčenko-Pastur Theorem, that the empirical covariance matrix of a high-dimensional white noise cannot be close to the identity matrix when the sample size $n$ and the number of variables $p$ have the same order (more strongly, that the average distance from its eigenvalues to $1$ has order $\sqrt{p/n}$). 
 \beg{propo}\la{flatMPcase} Let $M=\lf[x_{k,t}\ri]_{\substack{1\le k\le p\\ 1\le t\le n}}$ be a $p\ti n$ matrix of independent, centered, $L^4$ real random variables \st for all $k, t$, $\E[x_{k,t}^2]=1.$ Then we have $$\E\lf[\Frobnorm{\ff{n}MM^\trans -\bI_p}^2\ri]\le \f{p}n \lf( p+\max_{k,t}\E[x_{k, t}^4]-2\ri),$$ with   equality if all $x_{k,t}$'s have the same fourth moment.
 \end{propo}
 \bpr Given $\E(MM^\trans)=\bI_p$, the LHT reduces to $\ff{n^2}\E[\Tr(MM^\trans MM^\trans)] -p.$
 We have \beq \E[\Tr(MM^\trans MM^\trans)]&=&\sum_{1\le k, l\le p}\sum_{1\le t, s\le n} \E[x_{k,t}x_{l,t}x_{l,s}x_{k,s}]\eeq 
 By independence, the terms in this sum corresponding to ($k\ne l$ and $s\ne t$) are null. 
 The terms corresponding to ($k=l$ and $s=t$) have a total contribution upper-bounded by $np\max_{k,t}\E[x_{k, t}^4]$ (with equality if the fourth moments coincide).  The terms corresponding to ($k\ne l$ and $s=t$) (resp. to ($k=l$ and $s\ne t$)) have a total contribution equal to $np(p-1)$ (resp. $pn(n-1)$). One then concludes easily.
 \epr 
 
The following identity, which can be proved by induction, is useful when dealing with the squared Frobenius norm of a sum.
 \blem For any $a_1, \ld, a_k\in \R$, \be\la{ineqsumsquares}(a_1+\cdots+a_k)^2\le 2^{k-1}(a_1^2+\cdots+a_k^2).\ee
 \elem 
 
 \subsection{Variance in a sparse graphical model}
 In the following lemma, whose proof is straightforward,  the exponent $2d-1$ (instead of the natural $2d$) is what has most  of our convergences working.
 \blem Let  $$\lf(\sfy_{\bi}\ri)_\ijinnsqt$$ be an array  of $L^2$ real random variables \st for all $\bi=(i_1, \ld, i_d)\in \nsqt$, $\bj=(j_1, \ld, j_d)\in \nsqt$,   $$\{i_1, \ld, i_d\} \cap\{j_1, \ld, j_d\} =\emptyset \implies \Cov(\sfy_{\bi}, \sfy_{\bj})=0.$$
 Then \be\la{keyupperbound2308252}\Var\lf(\sum_\ijinnsqt\sfy_{\bi}\ri)\le 2n^{2d-1}\max_\ijinnsqt\Var(\sfy_{\bi}).\ee
 \elem

\subsection{McDiarmid-Azuma-Hoeffding concentration
inequality for bounded differences}
The following   lemma can be found, for example, in \cite{McDiarmid} or \cite[Th. 6.2]{MR3185193}. 
\blem\la{McDiarmidLemma} Let $X_1, \ld, X_n$ be random variances taking values in respectively $\cS_1, \ld, \cS_n$ and $f: \cS_1\ti \cd\ti\cS_n\to\R$ measurable be \st for any $i\in [n]$, \be\la{defsiisquare}\si_i=\sup_{\substack{x, x'\in \cS_1\ti \cd\ti\cS_n\\ \forall j\ne i, x_j=x'_j}}|f(x')-f(x)|\ee is finite. Then, for any $r\ge 0$, $$\p\lf(|f(X_1, \ld, X_n)-\E f(X_1, \ld, X_n)|\ge r\ri)\le 2\me^{-2r^2/\si^2}$$ for $\si^2$ defined by $$\si^2=\si_1^2+\cd+\si_n^2.$$
\elem

\section*{Acknowledgments}
We thank Victor Le Coz, Matteo Smerlak and Jean-Philippe Bouchaud, who contributed to our research through fruitful discussions.

This research was primarily conducted within the Econophysics \& Complex Systems Research Chair, under the aegis of the Fondation du Risque, the Fondation de l'Ecole Polytechnique, the Ecole Polytechnique and Capital Fund Management.

\bibliographystyle{plain}
\bibliography{sample}

\end{document}